\documentclass{article}

\usepackage{amssymb,latexsym,amsmath}

\usepackage[pdftex]{graphicx}

\usepackage{hyperref}

\begin{document}

\pagestyle{plain}

\newtheorem{theorem}{Theorem}

\newtheorem{proposition}[theorem]{Proposition}

\newtheorem{lemma}[theorem]{Lemma}

\newtheorem{corollary}[theorem]{Corollary}

\newtheorem{definition}[theorem]{Definition}

\newtheorem{remark}[theorem]{Remark}

\newtheorem{exempl}{Example}[section]

\newenvironment{exemplu}{\begin{exempl}  \em}{\hfill $\square$

\end{exempl}}  \vspace{.5cm}

\renewcommand{\contentsname}{ }

\title{Geometric Ruzsa triangle inequality in metric spaces with dilations}

\author{Marius Buliga \\ 
\\
Institute of Mathematics, Romanian Academy \\
P.O. BOX 1-764, RO 014700\\
Bucure\c sti, Romania\\
{\footnotesize Marius.Buliga@imar.ro}}  \vspace{.5cm}

\date{This version: 20.09.2016}

\maketitle

\begin{abstract}
This note  contains a general  "geometric Ruzsa triangle inequality"  in  metric spaces with dilations, inspired by lemma 4.1 and further comments from  the article \href{http://arxiv.org/abs/1212.5056}{arXiv:1212.5056 [math.CO]} "On growth in an abstract plane" by Nick Gill, H. A. Helfgott, Misha Rudnev. 
\end{abstract}

\paragraph{Motivation.} In \cite{unu} lemma 4.1  is given a proof of the Ruzsa triangle inequality which intrigued me. Later on, at the end of the article the authors give a geometric Ruzsa inequality in a Desarguesian projective plane, based on similar ideas as the ones used in the proof of the Ruzsa triangle inequality. I shall adapt the same idea to the frame of metric spaces with dilations, where: 
\begin{enumerate}
\item[-]  there is no algebraic structure, like the one of a group, except the approximate operations  provided by the field of dilations, 
\item[-] nor any incidence structure, like in a projective space, for example. 
\end{enumerate}


\paragraph{First result.} Let $ X$ be a non-empty set and $ \Delta: X \times X \rightarrow X$ be an operation on $ X$ which has the following two properties:
\begin{enumerate}
\item[1.] for any $ a, b, c \in X$ we have $ \Delta(\Delta(a,b), \Delta(a,c)) = \Delta(b,c)$, 
\item[2.] for any $ a \in X$   the function $ z \mapsto \Delta(z,a)$ is injective. 
\end{enumerate}

We may use weaker hypotheses for $ \Delta$, namely:
\begin{enumerate}
\item[1.](weaker) there is a function $ F: X \times X \rightarrow X$ such that $ F(\Delta(a,b), \Delta(a,c)) = \Delta(b,c)$ for any $ a, b, c \in X$, 
\item[2.](weaker) there is a function $ G: X \times X \rightarrow X$ such that $ a \mapsto G(\Delta(a,b), b)$ is an injective function for any $ b \in X$.
\end{enumerate}

\begin{proposition}
 Let $ X$ be a non empty set endowed with an operation $ \Delta$ which satisfies 1. and 2. (or the weaker version of those). Then for any non empty sets $ A, B, C \subset X$ there is an injection

$$ i: \Delta(C,A) \times B \rightarrow \Delta(B,C) \times \Delta(B,A)$$

where we denote by $ \Delta(A,B) = \left\{ \Delta(a,b) \mid a \in A, b \in B \right\}$

In particular, if $ A, B, C$ are finite sets, we have the Rusza triangle inequality

$$ \mid \Delta(C,A) \mid \mid B \mid \leq \mid \Delta(B,C) \mid \mid \Delta(B,A) \mid$$

where $ \mid A \mid$ denotes the cardinality of the finite set $ A$.
\label{prop1}
\end{proposition}

I shall give the proof for hypotheses 1. and 2. because the proof is the same for the weaker hypotheses. Also, this is basically the same proof as the one of the mentioned  lemma 4.1.  The proof of the Ruzsa inequality corresponds to the choice $ \Delta(a,b) = -a + b$, where $ (X,+)$ is a group (no need to be abelian). The proof  of the geometric Ruzsa inequality corresponds to the choice $ \Delta(a,b) = [b,a]$, with the notations from the article, with the observation that this function $ \Delta$ satisfies  the weaker versions of 1. and 2.

\paragraph{Proof.}  We can choose functions $ f: \Delta(C,A) \rightarrow C$ and $ g: \Delta(C,A) \rightarrow A$ such that for any $ x \in \Delta(C,A)$ we have $ x = \Delta(f(x),g(x))$. With the help of these functions let

$$ i(x,b) = (\Delta(b,f(x)), \Delta(b, g(x))) \quad \quad .$$
We want to prove that $ i$ is injective. Let $ (c,d) = i(x,b) = i(x',b')$. Then, by 1. we have $ x = x' = \Delta(c,d)$.  This gives an unique $ e = f(x) = f(x')$. Now we know that $ \Delta(b, e) = \Delta(b,f(x)) = c = \Delta(b', f(x')) = \Delta(b', e)$. By 2. we get that $ b = b'$     qed.


\paragraph{Metric spaces with dilations} These spaces appeared first  under the name "dilatation structures"  in the article \cite{buligadil1}.   In particular regular sub-riemannian manifolds, riemannian manifolds and Lie groups with a left invariant distance induced by a completely non-integrable (i.e. generating) distribution are examples of such spaces. For the most advanced presentation see  the course notes \cite{doi}. 

A metric space with dilations is a metric space $(X,d)$ endowed with locally defined groups of dilations: 
$$\delta_{\varepsilon}^{a}: \, dom \, \delta_{\varepsilon}^{a} \, \subset X \rightarrow X$$
for any point $a \in X$, where $\varepsilon \in \Gamma$ is an element of a topological commutative group $\Gamma$ with an absolute, called the group of scales. In many applications we take $\Gamma = (0,\infty)$ with multiplication and $0$ as absolute. In this article we use the same group, so $\varepsilon > 0$ is a generic element of the group $\Gamma$. 

For a generic scale $\varepsilon$ the various compositions of dilations can be described as binary trees with white and black nodes and with leaves decorated by the elements of the space $X$. Indeed, dilations appear as white or black nodes, according to the convention from the following figure. 

\vspace{.5cm}

\centerline{\includegraphics[width=  50mm]{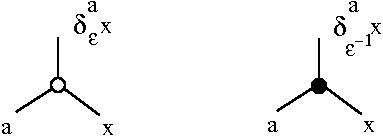}}

\vspace{.5cm}

The dilations indexed by the scale $1$ (i.e. the neutral element of the group $\Gamma$) are the identity. The dilations based at the same point compose as usual: 
$$\delta_{\varepsilon}^{a} \, \delta_{\mu}^{a} x \, = \, \delta_{\varepsilon \mu}^{a} x$$
whenever both members of the equality make sense. In the language of binary trees, for $\displaystyle \mu = \varepsilon^{-1}$ this transforms into a graph rewrite: 

\vspace{.5cm}

\centerline{\includegraphics[width=  50mm]{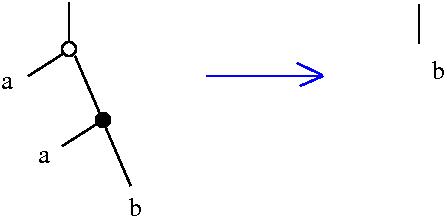}}

\vspace{.5cm}

 The base point of a dilation is a fixed point, i.e. for any $a \in X$  we have 
$$\delta_{\varepsilon}^{a} a \, = \, a$$
and any scale $\varepsilon > 0$. (There is another graph rewrite which corresponds to this, but it is not needed in this article.) 

In  a metric space with dilations  $ (X, d, \delta)$  we have the  approximate difference function, in graphical notation 

\vspace{.5cm}

\centerline{\includegraphics[width=  25mm]{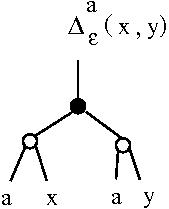}}

\vspace{.5cm}

$$ \Delta^{e}_{\varepsilon} (a,b) \, = \, \delta_{\varepsilon^{-1}}^{\delta_{\varepsilon}^{e} a} \, \delta_{\varepsilon}^{e} b$$ 
based at $ e \in X$ and applied to a pair of close points $ a, b \in X$. 

By the definition of a dilation structure, this function has the property that $\displaystyle (e,a,b) \mapsto \Delta^{e}_{\varepsilon}(a,b)$ converges uniformly (on compact sets) to $\displaystyle \Delta^{e}(a,b)$ as $ \varepsilon$ goes to $ 0$. 

An important result is that there exists a local group operation with $e$ as neutral element such that $ \Delta^{e}(a,b) = -a+b$, therefore the function $ \Delta^{e}$ satisfies 1. and 2. The local group operation is "conical", in the sense that the group of dilations based at $e$ are all morphisms. Typical examples of continuous conical groups are nilpotent groups. 

\begin{lemma}
The approximate difference operation $ \Delta^{e}_{\varepsilon}$  satisfies the following version of property 1.:
\begin{enumerate}
\item[1.](approximate) for any $ e, a, b, c \in X$ which are sufficiently close and for any $ \varepsilon \in (0,1)$ we have, with the notation $ a(\varepsilon) = \delta^{e}_{\varepsilon} a$,  the relation

$$ \Delta^{\delta^{e}_{\varepsilon} a}_{\varepsilon}(\Delta^{e}_{\varepsilon}(a,b), \Delta^{e}_{\varepsilon}(a,c)) = \Delta^{e}_{\varepsilon}(b,c) \quad \quad .$$
\end{enumerate}
\end{lemma}

\paragraph{Proof.} Apply the graph rewrite in the marked places from the following figure. 

\vspace{.5cm}

\centerline{\includegraphics[width=  100mm]{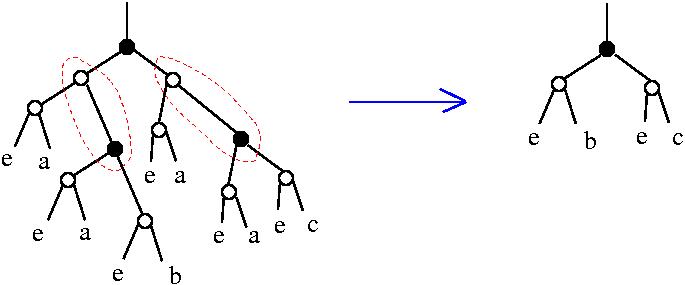}}

\vspace{.5cm}

We say that a set $ A \subset X$ is $ \varepsilon$ separated if for any $ x, y \in A$,  the inequality  $ d(x,y) < \varepsilon$  implies  $ x = y$.  Further I am going to write about sets which are close to a fixed, but arbitrary otherwise point $ e \in X$.

\begin{proposition}  In a metric space with dilations, let $ p > 0$ and  let $ A, B, C$ be finite sets of points included in a compact neighbourhood of $ e$, which are close to $ e \in X$, such that for any  $ \varepsilon \in (0,p)$  the sets $ B$ and $\displaystyle \Delta^{e}_{\varepsilon}(C,A)$ are $ \mu$ separated. Then for any $ \varepsilon \leq C(\mu)$ there is an injective function

$$ i^{e}_{\varepsilon}: \Delta^{e}_{\varepsilon}(C,A) \times B \rightarrow \Delta^{e}_{\varepsilon}(B,C) \times \Delta^{e}_{\varepsilon}(B,A) \quad \quad .$$
\label{prop2}
\end{proposition}

\paragraph{Proof.} As previously, we first choose  functions: 
$$f^{e}_{\varepsilon}: \Delta^{e}_{\varepsilon}(C,A) \rightarrow C \, , \, g^{e}_{\varepsilon}: \Delta^{e}_{\varepsilon}(C,A) \rightarrow A $$
such that 
$$\Delta^{e}_{\varepsilon}(f^{e}_{\varepsilon}(x), g^{e}_{\varepsilon}(x)) \, = \, x $$
for any $\displaystyle x \in \Delta^{e}_{\varepsilon}(C,A)$. We extend these functions by: 
$$f^{\delta^{e}_{\varepsilon} b}_{\varepsilon}(x) \, = \, \Delta^{e}_{\varepsilon}(b, f^{e}_{\varepsilon}(x)) $$
$$g^{\delta^{e}_{\varepsilon} b}_{\varepsilon}(x) \, = \, \Delta^{e}_{\varepsilon}(b, g^{e}_{\varepsilon}(x)) $$ 
for any $b \in B$ and $\displaystyle x \in \Delta^{e}_{\varepsilon}(C,A)$. With these functions we define: 
$$i^{e}_{\varepsilon}(x,b) \, = \, \left(f^{\delta^{e}_{\varepsilon} b}_{\varepsilon}(x) , g^{\delta^{e}_{\varepsilon} b}_{\varepsilon}(x) \right) \, \in \Delta^{e}_{\varepsilon}(B,C) \times \Delta^{e}_{\varepsilon}(B,A)$$
and we want to prove that this function is injective. 

Remark that by the property 1.(approximate) we have 
$$\Delta^{\delta_{\varepsilon}^{e} b}_{\varepsilon}(i^{e}_{\varepsilon}(x,b)) \, = \, x $$
therefore for any $(x,b)$ in the domain of $\displaystyle i^{e}_{\varepsilon}$ we get the estimate: 
$$d\left(x, \Delta^{e}_{\varepsilon}(i^{e}_{\varepsilon}(x,b)) \right) \, < \, \frac{\mu}{2}$$
provided that $\varepsilon \leq C(\mu)$. If $\displaystyle i^{e}_{\varepsilon}(x,b) \, = \, i^{e}_{\varepsilon}(x',b')$ then by the previous estimate 
$$d(x,x') \, < \, \mu$$
therefore $x = x'$ by the $\mu$ separation property of the set $\displaystyle \Delta^{e}_{\varepsilon}(C,A)$. We are left with the proof of $b = b'$. 

When $x = x'$ we have $\displaystyle \Delta^{e}_{\varepsilon}(b, z) = \Delta^{e}_{\varepsilon}(b',z)$ for $\displaystyle z  = f^{e}_{\varepsilon}(x)$. Again by the uniform convergence of the approximate difference $\displaystyle \Delta^{e}_{\varepsilon}(b,z)$ to a difference operation $\displaystyle \Delta^{e}(b,z)$ in a group, it follows that for a well chosen $C(\mu)$ we shall have $d(b,b') < \mu$ if $\varepsilon \leq C(\mu)$. By the $\mu$ separation property of the set $B$ we obtain $b = b'$.   qed.


\paragraph{Acknowledgements.} This article was supported by a grant of the Romanian National Authority for Scientific Research, project PN-II-ID-PCE-2011-3-0383.

\end{document}